\documentstyle[11pt,amssymb]{article}

\begin{document}

\begin{flushleft}
\Large{\bf Iterated integrals and unipotent periods on families of marked elliptic curves} 
\end{flushleft}

\begin{flushleft} 
\large{\bf Takashi Ichikawa} 
\end{flushleft}

\begin{flushleft} 
{\small Department of Mathematics, Faculty of Science and Engineering, Saga University, 
Saga 840-8502, Japan. E-mail: ichikawn@cc.saga-u.ac.jp} 
\end{flushleft} 

\noindent
{\bf Abstract:} \ 
We describe iterated integrals as unipotent periods on families of 
marked elliptic curves in terms of multiple zeta values and 
elliptic multiple zeta values. 
\vspace{2ex}

\vspace{2ex}

\vspace{2ex}

\noindent
{\bf 1. Introduction} 
\vspace{2ex} 

\noindent 
Recently, iterated integrals on families of marked elliptic curves were well studied 
in the study of Feynman integrals (for example, see \cite{W} and references therein) 
and of (nonabelian) unipotent extensions of elliptic period integrals 
(cf. \cite{CEE, E, H, L}). 
Considering results of of Brown \cite{Br} and Banks-Panzer-Pym \cite{BPP} 
in the case of marked Riemann surfaces, 
it is natural to pose a conjecture that these integrals are explicitly described by 
multiple zeta values and elliptic multiple zeta values 
which are expressed as iterated integrals of $dz/z, dz/(1-z)$ 
and of Eisenstein series respectively. 

The aim of this paper is to make a contribution to this conjecture. 
More precisely, 
for maximally degenerating families of marked elliptic curves over ${\mathbb C}$, 
we describe these unipotent periods as computable perturbative expansion 
in canonical deformation parameters whose coefficients are generated over ${\mathbb Q}$ by multiple zeta values and elliptic multiple zeta values. 
In order to obtain this description, 
we use formal geometry over ${\mathbb Z}$ which provides arithmetic 
gluing processes of marked Riemann spheres and elliptic curves 
using Grothendieck's (formal) existence theorem \cite[Section 5]{GD}. 
Note that this result be extended for curves of general genus 
by the theory of the universal Mumford curve given in \cite{I3}. 

In this paper, we first consider the universal marked Tate curve 
as an arithmetic universal family of marked elliptic curves, 
and construct the universal KZB connection on this curve 
by gluing the KZ (Knizhnik-Zamolodchikov) connection and its elliptic extension 
which are called the KZB (Knizhnik-Zamolodchikov-Bernard) connection. 
Therefore, we obtain {\it universal unipotent periods} which give isomorphisms 
between the Betti and de Rham tannakian fundamental groups 
for degenerating families of marked elliptic curves over ${\mathbb C}$. 
Then by studying the behavior of these periods under fusing moves on marked Tate curves, 
we show the main result stated above on the universal unipotent periods 
which is an extension of results of \cite{Br, BPP} in the case of genus $0$ 
to that of genus $1$. 
\vspace{4ex}

%{\it Acknowledgments.} 

\noindent
{\bf 2. Universal marked Tate curve} 
\vspace{2ex} 

\noindent
{\it 2.1. Marked Tate curve.} 
A (marked) curve is called {\it degenerate} if it is a stable (marked) curve and 
the normalization of its irreducible components are all projective (marked) lines. 
Then the dual graph $\Delta = (V, E, T)$ of a stable marked curve is a collection  
of 3 finite sets $V$ of vertices, $E$ of edges, $T$ of tails 
and 2 boundary maps 
$$
b : T \rightarrow V, 
\ \ b : E \longrightarrow \left( V \cup \{ \mbox{unordered pairs of elements of $V$} \} \right) 
$$
such that the geometric realization of $\Delta$ is connected 
and that $\Delta$ is {\it stable}, namely its each vertex has at least $3$ branches. 
Denote by $\sharp X$ the number of elements of a finite set $X$, 
and call a (connected) stable graph $\Delta = (V, E, T)$ {\it of $(g, n)$-type} 
if ${\rm rank}_{\mathbb Z} H_{1}(\Delta, {\mathbb Z}) = g$, $\sharp T = n$. 
Then under fixing a bijection $\nu : T \stackrel{\sim}{\rightarrow} \{ 1, ... , n \}$, 
which we call a numbering of $T$, 
$\Delta = (V, E, T)$ becomes the dual graph of a degenerate $n$-marked curve of genus $g$ 
such that each tail $h \in T$ corresponds to the $\nu(h)$th marked point. 
If $\Delta$ is trivalent, i.e. any vertex of $\Delta$ has just $3$ branches, 
then a degenerate $\sharp T$-marked curve with dual graph $\Delta$ 
is maximally degenerate. 
An {\it orientation} of a stable graph $\Delta = (V, E, T)$ means 
giving an orientation of each $e \in E$. 
Under an orientation of $\Delta$, 
denote by $\pm E = \{ e, -e \ | \ e \in E \}$ the set of oriented edges, 
and by $v_{h}$ the terminal vertex of $h \in \pm E$ (resp. the boundary vertex of $h \in T)$. 
For each $h \in \pm E$, 
denote by let $| h | \in E$ be the edge $h$ without orientation. 

Let $\Delta = (V, E, T)$ be a stable graph. 
Fix an orientation of $\Delta$, 
and take a subset ${\cal E}_{\rm f}$ of $\pm E \cup T$ 
whose complement ${\cal E}_{\infty}$ satisfies the condition that 
$$
\pm E \cap {\cal E}_{\infty} \cap \{ -h \ | \ h \in {\cal E}_{\infty} \} 
\ = \ 
\emptyset, 
$$ 
and that $v_{h} \neq v_{h'}$ for any distinct $h, h' \in {\cal E}_{\infty}$. 
We attach variables $x_{h}$ for $h \in {\cal E}_{\rm f}$ and $y_{e} = y_{-e}$ for $e \in E$ 
which are called moduli parameters and deformation parameters 
associated with $\Delta$. 
Let $R_{\Delta}$ be the ${\mathbb Z}$-algebra generated by $x_{h}$ 
$(h \in {\cal E}_{\rm f})$, 
$1/(x_{e} - x_{-e})$ $(e, -e \in {\cal E}_{\rm f})$ and $1/(x_{h} - x_{h'})$ 
$(h, h' \in {\cal E}_{\rm f}$ with $h \neq h'$ and $v_{h} = v_{h'})$, 
and let 
$$ 
A_{\Delta} \ = \ R_{\Delta} [[y_{e} \, (e \in E)]], \ \ 
B_{\Delta} \ = \ A_{\Delta} \left[ \prod_{e \in E} y_{e}^{-1} \right]. 
$$ 
For $h \in \pm E$, put 
\begin{eqnarray*}
\phi_{h} 
& = & 
\left( \begin{array}{cc} x_{h} & x_{-h} \\ 1 & 1 \end{array} \right) 
\left( \begin{array}{cc} 1 & 0 \\ 0 & y_{h} \end{array} \right) 
\left( \begin{array}{cc} x_{h} & x_{-h} \\ 1 & 1 \end{array} \right)^{-1} 
\\ 
& = & 
\frac{1}{x_{h} - x_{-h}} 
\left\{ \left( \begin{array}{cc} x_{h} & - x_{h} x_{-h} \\ 1 & - x_{-h} \end{array} \right) - 
\left( \begin{array}{cc} x_{-h} & - x_{h} x_{-h} \\ 1 & - x_{h} \end{array} \right) y_{h} \right\}, 
\end{eqnarray*} 
where $x_{h}$ (resp. $x_{-h})$ means $\infty$ 
if $h$ (resp. $-h)$ belongs to ${\cal E}_{\infty}$. 
This gives an element of $PGL_{2}(B_{\Delta}) = GL_{2}(B_{\Delta})/B_{\Delta}^{\times}$ 
denoted by the same symbol which satisfies
$$
\frac{\phi_{h}(z) - x_{h}}{z - x_{h}} 
\ = \ 
y_{h} \frac{\phi_{h}(z) - x_{-h}}{z - x_{-h}} 
\ \ (z \in {\mathbb P}^{1}), 
$$
where $PGL_{2}$ acts on ${\mathbb P}^{1}$ by linear fractional transformation. 
For any path $\rho = h(1) \cdot h(2) \cdots h(l)$ in $\Delta$ 
which is reduced in the sense that $h(i) \neq - h(i+1)$, 
one can associate an element $\rho^{*}$ of $PGL_{2}(B_{\Delta})$ 
having reduced expression $\phi_{h(l)} \phi_{h(l-1)} \cdots \phi_{h(1)}$. 
Then it is shown in \cite[Proposition 1.3]{I1} and its proof that 
there exist elements $\alpha, \alpha' \in B_{\Delta}$ and 
$\beta \in \left( \prod_{i=1}^{l} y_{h(i)} \right) \cdot (A_{\Delta})^{\times}$ 
such that 
$$
\frac{\rho^{*}(z) - \alpha}{z - \alpha} = 
\beta \frac{\rho^{*}(z) - \alpha'}{z - \alpha'} \ \ (z \in {\mathbb P}^{1}). 
$$ 
We call $\alpha, \alpha'$ the {\it attractive, repulsive fixed points} of $\rho^{*}$ 
respectively, 
and call $\beta$ the {\it multiplier} of $\rho^{*}$. 
\vspace{2ex}

\noindent
{\bf Theorem 2.1.} 
\begin{it} 
Let $\Delta = (V, E, T)$ be a stable graph of $(1, n)$-type. 
Then there exists a stable marked curve ${\cal E}_{\Delta}$ of genus $1$, 
called a $n$-marked Tate curve, over $A_{\Delta}$ 
which satisfies the following properties: 

\begin{itemize}

\item[\rm (P1)] 
The closed fiber ${\cal E}_{\Delta} \otimes_{A_{\Delta}} R_{\Delta}$ of 
${\cal E}_{\Delta}$ 
obtained by substituting $y_{e} = 0$ $(e \in E)$ 
becomes the degenerate marked curve over $R_{\Delta}$ with dual graph $\Delta$ which is 
obtained from the collection of $P_{v} := {\mathbb P}^{1}_{R_{\Delta}}$ $(v \in V)$ 
by identifying the points $x_{e} \in P_{v_{e}}$ and $x_{-e} \in P_{v_{-e}}$ ($e \in E$), 
where $x_{h}$ denotes $\infty$ if $h \in {\cal E}_{\infty}$. 

\item[\rm (P2)] 
${\cal E}_{\Delta}$ gives rise to a universal deformation of 
${\cal E}_{\Delta} \otimes_{A_{\Delta}} R_{\Delta}$ by the deformation parameters 
$y_{e}$ $(e \in E)$. 

\item[\rm (P3)] 
Taking $x_{h}$ $(h \in {\cal E}_{\rm f})$ as complex numbers such that 
$x_{e} \neq x_{-e}$ and that $x_{h} \neq x_{h'}$ if $h \neq h'$ and $v_{h} = v_{h'}$, 
and $y_{e}$ $(e \in E)$ as sufficiently small nonzero complex numbers, 
${\cal E}_{\Delta}$ gives rise to a family of $n$-marked elliptic curves over ${\mathbb C}$. 

\end{itemize}
\end{it} 

\noindent
{\it Proof.} 
When $n =1$, $\sharp V = \sharp E = 1$ and 
$$
\{ x_{h} \ | \ h \in \pm E \cup T \} = \{ 0, 1, \infty \}, 
$$
$R_{\Delta} = {\mathbb Z}$ and $A_{\Delta} = {\mathbb Z}[[ y_{e} ]]$, 
where $E = \{ e \}$. 
Then it is shown by Tate (cf. \cite{T, Si}) that ${\cal E}_{\Delta}$ is given as 
the Tate (elliptic) curve over ${\mathbb Z}[[ y_{e} ]]$ 
formally represented as ${\mathbb G}_{m}/ \langle y_{e} \rangle$ which becomes 
the quotient space of ${\mathbb C}^{\times} = {\mathbb C} - \{ 0 \}$ 
under the multiplication by a complex number $y_{e}$ with $|y_{e}| < 1$. 
Therefore, if there is a loop in $E$ which we denote by $l$, 
then there exists a formal  scheme $\widetilde{\cal E}_{\Delta}$ over $A_{\Delta}$ 
which is obtained from the Tate curve ${\mathbb G}_{m}/ \langle y_{l} \rangle$ 
and $P_{v} = {\mathbb P}^{1}$ $(v \in V)$ with marked points $x_{t}$ 
(\mbox{$t \in T$ with $v_{t} = v$}) by gluing them with $\phi_{e}$ $(e \in E - \{ l \})$. 
Since $\widetilde{\cal E}_{\Delta} \otimes_{A_{\Delta}} R_{\Delta}$ has an ample divisor, 
by Grothendieck's existence theorem, 
we have a required marked curve ${\cal E}_{\Delta}$ 
as the algebraization of $\widetilde{\cal E}_{\Delta}$, 
namely $\widetilde{\cal E}_{\Delta}$ is obtained as the formal completion of ${\cal E}_{\Delta}$ 
along its closed fiber. 

When there is no loop in $E$, 
as is shown by Ihara-Nakamura \cite{IhN} in the case of general genus, 
by Grothendieck's existence theorem, 
${\cal E}_{\Delta}$ is given as the algebraization of the formal scheme 
over $A_{\Delta}$ which is obtained from $P_{v} = {\mathbb P}^{1}$ $(v \in V)$ 
with marked points $x_{t}$  (\mbox{$t \in T$ with $v_{t} = v$}) by gluing them 
with $\phi_{e}$ $(e \in E)$. 
\ $\square$ 
\vspace{2ex}

\noindent
{\it 2.2. Comparison of parameters.}
Let $\Delta_{1} = (V_{1}, E_{1}, T_{1})$ be a stable graph of $(1, n)$-type 
which is not trivalent. 
Then there exists a vertex $v_{0} \in V_{1}$ which has at least $4$ branches. 
Take two elements $h_{1}, h_{2}$ of $\pm E_{1} \cup T_{1}$ such that 
$h_{1} \neq h_{2}$ and $v_{h_{1}} = v_{h_{2}} = v_{0}$, 
and let $\Delta_{2} = (V_{2}, E_{2}, T_{2})$ be a stable graph obtained from $\Delta_{1}$ 
by replacing $v_{0}$ with an oriented (nonloop) edge $h_{0}$ such that 
$v_{h_{1}} = v_{h_{2}} = v_{h_{0}}$ and that $v_{h} = v_{-h_{0}}$ 
for any $h \in \pm E_{1} \cup T_{1} - \{ h_{1}, h_{2} \}$ with $v_{h} = v_{0}$. 
Put $e_{i} = |h_{i}|$ for $i = 0, 1, 2$. 
Then we have the following identifications: 
$$
V_{1} = V_{2} - \{ v_{-h_{0}} \} \ (\mbox{in which $v_{0} = v_{h_{0}}$}), \ 
E_{1} = E_{2} - \{ e_{0} \}, \ T_{1} = T_{2}, 
$$ 
and denote by $s_{e}$ the deformation parameters corresponding to 
$e \in E_{2} = E_{1} \cup \{ e_{0} \}$. 
\vspace{2ex} 

\noindent
{\bf Theorem 2.2.} 
\begin{it} 

{\rm (1)} 
The marked Tate curves ${\cal E}_{\Delta_{1}}, {\cal E}_{\Delta_{2}}$ 
associated with $\Delta_{1}, \Delta_{2}$ respectively 
are isomorphic over $R_{\Delta_{2}} \left[ s_{e_{0}}^{-1} \right] [[ s_{e} \, (e \in E_{2}) ]]$, 
where 
$$
\frac{x_{h_{1}} - x_{h_{2}}}{s_{e_{0}}}, \ \ 
\frac{y_{e_{i}}}{s_{e_{0}} s_{e_{i}}} \ (\mbox{$i = 1, 2$ with $h_{i} \not\in T_{1}$}), \ \
\frac{y_{e}}{s_{e}} \ (e \in E_{1} - \{ e_{1}, e_{2} \}) 
$$
belong to $(A_{\Delta_{2}})^{\times}$ if $h_{1} \neq - h_{2}$, 
and 
$$
\frac{x_{h_{1}} - x_{h_{2}}}{s_{e_{0}}}, \ \ 
\frac{y_{e}}{s_{e}} \ (e \in E_{1}) 
$$
belong to $(A_{\Delta_{2}})^{\times}$ if $h_{1} = - h_{2}$. 

{\rm (2)} 
The assertion (1) holds in the category of complex geometry 
when $x_{h_{1}} - x_{h_{2}}, y_{e}$ and $s_{e}$ are taken to be 
sufficiently small complex numbers. 

\end{it}
\vspace{2ex}

\noindent
{\it Proof.} 
We prove the assertion (1). 
By \cite[Lemma 1.2]{I1}, 
$\phi_{-h_{0}}(t_{h_{1}}) - \phi_{-h_{0}}(t_{h_{2}})$ belongs to 
$s_{e_{0}} \cdot (A_{\Delta_{2}})^{\times}$, 
and hence 
${\cal E}_{\Delta_{2}} \otimes_{A_{\Delta_{2}}} 
R_{\Delta_{2}} \left[ s_{e_{0}}^{-1} \right] [[ s_{e} \ (e \in E_{2} )]]$ 
gives a universal deformation of a universal degenerate curve with dual graph $\Delta_{1}$. 
Then by the universality of marked Tate curves, 
there exists an injective homomorphism 
$A_{\Delta_{1}} \hookrightarrow R_{\Delta_{2}} \left[ s_{e_{0}}^{-1} \right] [[ s_{e} ]]$
which gives rise to an isomorphism ${\cal E}_{\Delta_{1}} \cong {\cal E}_{\Delta_{2}}$. 
Denote by $t_{h}$ the moduli parameters of ${\cal E}_{\Delta_{2}}$ corresponding to 
$h \in \pm E_{2} \cup T_{2}$.  
Then under this homomorphism, 
\begin{eqnarray*}
\lefteqn{
\left( P_{v_{-h_{0}}}; \phi_{-h_{0}}(t_{h_{1}}), \phi_{-h_{0}}(t_{h_{2}}), 
t_{h} \ (v_{h} = v_{-h_{0}}, h \neq h_{0}) \right) 
} 
\\ 
& \cong & 
\left( P_{v_{0}}; x_{h_{1}}, x_{h_{2}}, x_{h} \ (v_{h} = v_{0}, h \neq h_{1}, h_{2}) \right), 
\end{eqnarray*} 
and hence $x_{h_{1}} - x_{h_{2}} \in s_{e_{0}} \cdot (A_{\Delta_{2}})^{\times}$. 
Furthermore, when $h_{1} \neq -h_{2}$, 
the deformation parameters of 
$$
\left( P_{v_{-h_{0}}}; \phi_{-h_{0}}(t_{h_{1}}), \phi_{-h_{0}}(t_{h_{2}}), 
t_{h} \ (v_{h} = v_{-h_{0}}, h \neq h_{0}) \right) 
$$
corresponding to $h_{i} \cdot (-h_{0})$ $(i = 1, 2)$ are $y_{h_{i}}$, 
and hence by \cite[Proposition 1.3]{I1}, 
$y_{h_{i}} \in \left( s_{h_{0}} \cdot s_{h_{i}} \right) \cdot  (A_{\Delta_{2}})^{\times}$. 
When $h_{1} = - h_{2}$, 
the deformation parameters of 
$$
\left( P_{v_{-h_{0}}}; \phi_{-h_{0}}(t_{h_{1}}), \phi_{-h_{0}}(t_{h_{2}}), 
t_{h} \ (v_{h} = v_{-h_{0}}, h \neq h_{0}) \right) 
$$
corresponding to $h_{0} \cdot h_{1} \cdot (-h_{0})$ is $y_{h_{1}}$, 
and hence $y_{h_{1}} \in s_{h_{1}} \cdot  (A_{\Delta_{2}})^{\times}$. 
The assertion (2) follows from Theorem 2.1 (P3). 
\ $\square$ 
\vspace{2ex}

\noindent
{\it 2.3. Universal marked Tate curve.} 
For a positive integer $n$, 
denote by $\overline{\cal M}_{1,n}$ the moduli stack over ${\mathbb Z}$ 
of stable $n$-marked curves of genus $1$, 
and by ${\cal M}_{1,n}$ its substack classifying $n$-marked proper smooth curves of genus $1$ 
\cite{DM, KM, K}. 
Then by definition, there exits the universal stable $n$-marked curve 
${\cal E}_{1,n}$ of genus $1$ over $\overline{\cal M}_{1,n}$. 
\vspace{2ex}

\noindent
{\bf Theorem 2.3.} 
\begin{it} 
Let ${\cal D}_{1,n}$ be the closed substack of $\overline{\cal M}_{1,n}$ 
consisting of $n$-marked degenerate curves of genus $1$, 
and denote by ${\cal N}_{1,n}$ the formal completion of $\overline{\cal M}_{1,n}$ 
along ${\cal D}_{1,n}$. 
Then there exists an algebraization ${\cal N}^{\rm alg}_{1,n}$ of ${\cal N}_{1,n}$, 
namely ${\cal N}^{\rm alg}_{1,n}$ is a scheme containing ${\cal D}_{1,n}$ 
as its closed subset  such that ${\cal N}_{1,n}$ is the formal completion 
of ${\cal N}^{\rm alg}_{1,n}$ along ${\cal D}_{1,n}$, 
and the fiber of ${\cal E}_{1,n}$ over ${\cal N}_{1,n}$ 
gives a universal family of $n$-marked Tate curves.  
\end{it} 
\vspace{2ex}

\noindent
{\it Proof.} 
Let $\Delta = (V, E, T)$ be a stable graph of $(1, n)$-type, 
and take a system of coordinates on $P_{v} = {\mathbb P}^{1}_{R_{\Delta}}$ $(v \in V)$ 
such that $x_{h} = \infty$ $(h \in {\cal E}_{\infty})$ and that 
$\{ 0, 1 \} \subset P_{v}$ is contained in the set of points 
given by $x_{h}$ $(h \in {\cal E}_{\rm f}$ with $v_{h} = v)$. 
Under this system of coordinates, 
one has the marked Tate curve ${\cal E}_{\Delta}$ whose closed fiber 
${\cal E}_{\Delta} \otimes_{A_{\Delta}} R_{\Delta}$ gives a family of 
degenerate curves over the open subspace of 
$$
S = \left\{ \left. \left( p_{h} \in P_{v_{h}} \right)_{h \in \pm E \cup T} \ \right| \ 
p_{h} \neq p_{h'} \, (h \neq h', v_{h} = v_{h'}) \right\} 
$$ 
defined as $p_{e} \neq p_{-e}$ $(e \in E)$. 
Therefore, by taking another system of coordinates and 
comparing the associated marked Tate curves with 
the original ${\cal E}_{\Delta}$ in a similar way to the proof of Theorem 2.2, 
${\cal E}_{\Delta}$ can be extended over the algebraization 
of the formal completion of $S$ in $\overline{\cal M}_{1,n}$. 
Then by Theorem 2.1, 
marked Tate curves ${\cal E}_{\Delta}$ for various stable graphs $\Delta$ 
of $(1, n)$-type are glued over ${\cal N}_{1,n}$, 
and hence the assertion hods. 
\ $\square$  
\vspace{2ex} 

\noindent
{\it Definition 2.3.} 
We call this universal family of $n$-marked Tate curves 
the {\it universal $n$-marked Tate curve}. 
\vspace{4ex}

\noindent
{\bf 3. Unipotent periods of marked elliptic curves}
\vspace{2ex}

\noindent
{\it 3.1. Unipotent periods of Riemann surfaces.}
Fix nonnegative integers $g, n$ such that $2g - 2 + n > 0$. 
Let $R^{\circ}$ be a (not necessarily compact) Riemann surface 
obtained from a compact Riemann surface $R$ of genus $g$ by removing $n$ points. 
Then the category of unipotent local systems over ${\mathbb C}$ on $R^{\circ}$ 
is equivalent to that of vector bundles with nilpotent flat connection on $R$ 
which have regular singularities at $R - R^{\circ}$ with nilpotent residues (cf. \cite{D1}). 
We consider the tannakian category of unipotent local systems over ${\mathbb Q}$ 
on $R^{\circ}$ with fiber functor obtained from taking the fiber over 
a (tangential) point $x \in R^{\circ}$. 
Then its tannakian fundamental group $\pi^{\rm Be}_{1} (R^{\circ}, x)$ is 
a profinite algebraic group over ${\mathbb Q}$, 
and there exists a canonical homomorphism from the topological fundamental group 
$\pi_{1} (R^{\circ}, x)$ into $\pi^{\rm Be}_{1} (R^{\circ}, x)$. 

For a subfield $K$ of ${\mathbb C}$, 
let $C^{\circ}$ be a smooth curve obtained from a proper smooth curve $C$ 
of genus $g$ over $K$ by removing $n$ rational points over $K$. 
We consider the tannakian category of vector bundles with nilpotent flat connection 
over $K$ on $C$ which have regular singularities at $C - C^{\circ}$ with nilpotent residues 
with fiber functor obtained from taking the fiber over 
a $K$-rational (tangential) point $x \in C^{\circ}$. 
Then its tannakian fundamental group $\pi^{\rm dR}_{1} (C^{\circ}, x)$ is 
a profinite algebraic group over $K$, 
and by the above categorical equivalence, 
there exists a canonical isomorphism 
$$
\pi^{\rm dR}_{1} (C^{\circ}, x) \otimes_{K} {\mathbb C} \cong 
\pi^{\rm Be}_{1} ((C^{\circ})^{\rm an}, x) \otimes_{\mathbb Q} {\mathbb C}, 
$$
where $(C^{\circ})^{\rm an}$ denotes the Riemann surface associated with $C^{\circ}$. 
We call this isomorphism the {\it unipotent period isomophism} which is described by 
monodromy representations of $\pi^{\rm Be}_{1} ((C^{\circ})^{\rm an}, x)$. 
Similarly, for (tangential) points $x, y$ on $C^{\circ}$, 
we have the associated monodromy of a nilpotent flat connection. 
\vspace{2ex}

\noindent
{\it 3.2. KZ connection and connection matrix.} 
A {\it KZ connection} is defined on ${\mathbb P}^{1}$ is defined as a trivial bundle 
with flat connection which has regular singularities along $p_{i} \in {\mathbb P}^{1}$ 
with residue $X_{i}$ $(1 \leq i \leq m)$, 
where $X_{i}$ are symbols satisfying $\sum_{i=1}^{m} X_{i} = 0$. 
When $m = 3, p_{1} = 0, p_{2} = 1, p_{3} = \infty$, 
this connection is 
$$
df - f \left( X_{0} \frac{dz}{z} + X_{1} \frac{dz}{z-1} \right). 
$$
Then the associated monodromy from the tangential point $\vec{v}_{0} = d/dz$ at $0$ 
to $z$ with $|z| < 1$ becomes a noncommutative formal power series in $X_{0}, X_{1}$ 
whose coefficients are multiple polylogarithm functions 
$$
{\rm Li}_{k_{1},..., k_{l}}(z) = 
\int_{0}^{z} \frac{dz}{1-z} \underbrace{\frac{dz}{z} \cdots \frac{dz}{z}}_{k_{1}-1} 
\frac{dz}{1-z} \cdots \frac{dz}{1-z} \underbrace{\frac{dz}{z} \cdots \frac{dz}{z}}_{k_{l}-1} 
$$
of $z \in {\mathbb P}^{1} - \{ 0, 1, \infty \}$ regularized at $0$. 
Therefore, 
$$
{\rm Li}_{k_{1},..., k_{l}}(z) = 
\sum_{0 < n_{1} < \cdots < n_{l}} \frac{z^{n_{l}}}{n_{1}^{k_{1}} \cdots n_{l}^{k_{l}}} 
$$
are power series in $z$ over ${\mathbb Q}$. 

If $z$ is replaced with the tangential point $-\vec{v}_{1} = - d/dz$ at $1$, 
then the associated monodromy is called the {\it Drinfeld associator} \cite{D2, Dr} 
which is a noncommutative formal power series $\Phi (X_{0}, X_{1})$ in $X_{0}, X_{1}$ 
with coefficients expressed by multiple zeta values:  
$$
\zeta(k_{1},..., k_{l}) = 
\sum_{0 < n_{1} < \cdots < n_{l}} \frac{1}{n_{1}^{k_{1}} \cdots n_{l}^{k_{l}}} \ \ (k_{l} > 1). 
$$
The Drinfeld associator can be applied to calculating connection matrices 
for the KZ connection on the moduli space of a projective line with $4$ marked points 
$0, 1, \infty$ and $x$. 
When the associated connection has simple poles at $x = 0$ (resp. $x = 1$) 
with residue $X_{0}$ (resp. $X_{1}$), 
the connection matrix for the fusing move of $x \in (0, 1)$ 
between tangential points at $0$ and $1$ becomes $\Phi(X_{0}, X_{1})$. 
\vspace{2ex}

\noindent
{\it 3.3. KZB connection on the Tate curve.} 
We review results of Calaque-Enriquez-Etingof \cite{CEE}, Hain \cite{H} and Luo \cite{L} 
on the (universal) elliptic KZB connection following the formulation of \cite{H}. 
Let $q$ be a variable, 
and denote by $E_{q} = {\mathbb G}_{m}/ \langle q \rangle$ the Tate (elliptic) curve 
over ${\mathbb Z}[[q]]$. 
Then its closed fiber $E_{0} = E_{q} \otimes_{{\mathbb Z}[[q]]} {\mathbb Z}[[q]]/(q)$ 
minus the identity $1$ is identified with ${\mathbb P}^{1} - \{ 0, 1, \infty \}$. 
Denote by ${\mathbb Q} \left\langle\!\left\langle T, A \right\rangle\!\right\rangle$ 
of noncommutative formal power series over ${\mathbb Q}$ in the symbols $T, A$. 
Then under taking $\{ T, A \}$ as a de Rham framing of the first cohomology group of $E_{q}$, 
Hain \cite{H} constructs a vector bundle with flat connection, 
called the {\it elliptic KZB connection}, on $E_{q}$ whose fibers are identified with 
${\mathbb Q} \left\langle\!\left\langle T, A \right\rangle\!\right\rangle$ 
such that the associated flat connection has regular singularity at $1$ with residue $[T, A]$. 
Put $\vec{v}_{1} = d/dz$ at $1 \in E_{q} = {\mathbb G}_{m}/\langle q \rangle$. 
Then the monodromy of the elliptic KZB connection from $- \vec{v}_{1}$ to $\vec{v}_{1}$ 
in $E_{q}$ is called the Enriquez elliptic associator, 
and is represented as a noncommutative formal power series in $T, A$ 
whose coefficients are called {\it elliptic multiple zeta values}. 
The elliptic KZB connections is described in \cite[9.2]{H}, 
and it gives the above KZ connection on $E_{0} - \{ 1 \} = {\mathbb P}^{1} - \{ 0, 1, \infty \}$  by composing with the homomorphism 
$$
{\mathbb Q} \left\langle\!\left\langle X_{0}, X_{1}, X_{\infty} \right\rangle\!\right\rangle 
/ (X_{0} + X_{1} + X_{\infty}) 
\rightarrow 
{\mathbb Q} \left\langle\!\left\langle T, A \right\rangle\!\right\rangle 
$$
given by 
$$
X_{0} \mapsto \left( \frac{T}{e^{T} - 1} \right) \cdot A, \ 
X_{1} \mapsto [T, A], \ 
X_{\infty} = \left( \frac{T}{e^{-T} - 1} \right) \cdot A, 
$$
where $f(T, A) \cdot x = f \left( {\rm ad}_{T}, {\rm ad}_{A} \right)(x)$ 
for $f(T, A) \in {\mathbb Q} \left\langle\!\left\langle T, A \right\rangle\!\right\rangle$ 
(cf. \cite[Section 18]{H}). 
\vspace{2ex}

\noindent
{\it 3.4. KZB connection on marked Tate curves.} 
Let $\Delta_{0} = (V_{0}, E_{0}, T_{0})$ be a stable graph with orientation of $(1, n)$-type 
which consists of two vertices $v_{0}, v_{1}$, 
one (nonloop) edge $e$ and loop $l$ with orientation such that 
$v_{e} = v_{0}$, $v_{\pm l} = v_{-e} = v_{1}$. 
Then we call the associated marked Tate curve ${\cal E}_{\Delta_{0}}$ 
the {\it basic marked Tate curve}. 
We construct a vector bundle with flat connection $({\cal V}_{v_{0}}, {\cal F}_{v_{0}})$ 
on ${\cal E}_{\Delta_{0}}$ whose each fiber is identified with the ring 
$$
{\cal A}_{\Delta_{0}} = 
{\mathbb Q} \left\langle\!\left\langle X_{t}, T_{l}, A_{l} \right\rangle\!\right\rangle 
$$
of noncommutative formal power series over ${\mathbb Q}$ in the symbols 
$X_{t}$ $(t \in T_{0})$ and $T_{l}, A_{l}$ satisfying the condition 
$$
\sum_{t \in T} X_{t} - [T_{l}, A_{l}] = 0. 
$$
Put 
$$
X_{l} = \left( \frac{T_{l}}{e^{T_{l}} - 1} \right) \cdot A_{l}, \ 
X_{-l} = \left( \frac{T_{l}}{e^{-T_{l}} - 1} \right) \cdot A_{l}, 
$$
and $X_{e} = - X_{-e} = - \left[ T_{l}, A_{l} \right]$. 
Then for each vertex $v \in V_{0}$, 
the sum of $X_{h}$ for $h \in \pm E_{0} \cup T_{0}$ with $v_{h} = v$ is equal to $0$. 
\vspace{2ex}

\noindent
{\bf Theorem 3.1.} 
\begin{it} 
There exists a vector bundle with flat connection 
$\left( {\cal V}_{\Delta_{0}}, {\cal F}_{\Delta_{0}} \right)$ on ${\cal E}_{\Delta_{0}}$ 
which is constructed by gluing the KZ connection on $P_{v} = {\mathbb P}^{1}$ $(v \in V_{0})$
with regular singularities along $x_{h}$ with residue $X_{h}$ 
$(\mbox{$h \in \pm E_{0} \cup T_{0}$ such that $v_{h} = v_{0}$})$, 
and the elliptic KZB connections on the Tate curves associated with $l$. 
\end{it}
\vspace{2ex}

\noindent
{\it Proof.} 
Denote by $({\cal V}_{v_{0}}, {\cal F}_{v_{0}})$ the vector bundle with flat connection 
on $P_{v_{0}}$, 
where ${\cal V}_{v_{0}}$ is the trivial bundle with fiber ${\cal A}_{\Delta_{0}}$, 
and ${\cal F}_{v_{0}}$ has regular singularity along $x_{h}$ with residue $X_{h}$ 
for $h \in \pm E_{0} \cup T_{0}$ such that $v_{h} = v_{0}$. 
As reviewed in 3.3, 
there exists the elliptic KZB connection $({\cal V}_{l}, {\cal F}_{l})$ 
on the Tate elliptic curve $E_{q_{l}} = {\mathbb G}_{m}/ \langle q_{l} \rangle$ 
over ${\mathbb Z}[[ q_{l} ]]$, 
where fibers of ${\cal V}_{l}$ are identified with 
${\mathbb Q} \left\langle\!\left\langle T_{l}, A_{l} \right\rangle\!\right\rangle 
\hookrightarrow {\cal A}_{\Delta_{0}}$. 

Let $\xi_{e}$ (resp. $\xi_{-e}$) be tangential points over ${\mathbb Z}$ at 
$x_{e} \in P_{v_{0}} ={\mathbb P}^{1}$ (resp. the identity in $E_{q_{l}}$). 
Then the marked Tate curve ${\cal E}_{\Delta_{0}}$ can be obtained 
from $P_{v_{0}} \cup E_{q}$  defined by $\xi_{e} \cdot \xi_{-e} = y_{e}$, 
where $y_{e}$ are variables, 
and hence ${\cal E}_{\Delta_{0}}$ is proper over the ring $A_{\Delta_{0}}$ identified with 
$$
{\mathbb Z} \left[ x_{h}, \frac{1}{x_{h}- x_{h'}} \ 
(\mbox{$h, h' \in \{ e \} \cup T_{0}$ with $h \neq h'$}) \right] [[ q_{l}, y_{e} ]]. 
$$
Since 
$\displaystyle X_{e} \frac{d \xi_{e}}{\xi_{e}} = X_{-e} \frac{d \xi_{-e}}{\xi_{-e}}$, 
by gluing 
$({\cal V}_{v_{0}}, {\cal F}_{v_{0}})$ and $({\cal V}_{l_{i}},  {\cal F}_{l_{i}})$ $(1 \leq i \leq g)$, 
we obtain a vector bundle with flat connection on the formal completion 
$\widehat{\cal E}_{\Delta_{0}}$ of ${\cal E}_{\Delta_{0}}$ along its closed subscheme 
${\cal E}_{\Delta_{0}} \otimes_{A_{\Delta_{0}}} \left( A_{\Delta_{0}}/I \right)$, 
where $I$ is the ideal of $A_{\Delta_{0}}$ generated by $q_{l}, y_{e}$. 
Then by Grothendieck's existence theorem, 
there exists the associated vector bundle with flat connection on 
${\cal E}_{\Delta_{0}}$. 
\ $\square$ 
\vspace{2ex}

For each stable graph graph $\Delta = (V, E, T)$ of $(1, n)$-type with orientation, 
we attach symbols $X^{\Delta}_{h}$ $(h \in \pm E \cup T)$ 
by the following rules: 
\begin{itemize}

\item 
If $\Delta = \Delta_{0}$, 
then $X^{\Delta}_{h} = X_{h}$. 

\item 
For any $v \in V$, 
the sum of $X^{\Delta}_{h}$ $(h \in \pm E \cup T$ with $v_{h} = v)$ is $0$. 

\item 
Assume that $\Delta_{1} = (V_{1}, E_{1}, T_{1})$ and 
$\Delta_{2} = (V_{2}, E_{2}, T_{2})$ are given in 2.2. 
Then 
$$
X^{\Delta_{1}}_{h_{i}} = X^{\Delta_{2}}_{h_{0}} + X^{\Delta_{2}}_{h_{i}} \ (i = 1, 2), \ \ 
X^{\Delta_{1}}_{h} = X^{\Delta_{2}}_{h} \ (h \neq h_{1}, h_{2})
$$
if $h_{1} \neq -h_{2}$, and 
$$
X^{\Delta_{1}}_{h} = X^{\Delta_{2}}_{h} \ (h \neq h_{0}) 
$$
if $h_{1} = -h_{2}$.    

\end{itemize}

\noindent
{\bf Theorem 3.2.} 
\begin{it} 
The vector bundle with flat connection $({\cal V}_{\Delta_{0}}, {\cal F}_{\Delta_{0}})$ 
on ${\cal E}_{\Delta_{0}}$ can be analytically continued to a vector bundle with flat connection 
on the universal $n$-marked Tate curve given in Definition 2.3. 
This restriction to the closed fiber of ${\cal E}_{\Delta}$ 
for a stable graph $\Delta = (V, E, T)$ of $(1, n)$-type gives the KZ connection 
on $P_{v}$ $(v \in V)$ having regular singularities along $x_{h}$ with residue $X^{\Delta}_{h}$ 
for $h \in \pm E \cup T$ such that $v_{h} = v$. 
\end{it} 
\vspace{2ex}

\noindent
{\it Proof.} 
Let $\Delta_{1} = (V_{1}, E_{1}, T_{1}), 
\Delta_{2} = (V_{2}, E_{2}, T_{2})$ and $v_{0} \in V_{1}, h_{0} \in E_{2}$ be as in 2.2. 
Then $X^{\Delta_{2}}_{h_{0}} + X^{\Delta_{2}}_{-h_{0}} = 0$, 
and hence the KZ connection on $P_{v_{0}}$ is obtained from the KZ connections 
on $P_{v_{h_{0}}}, P_{v_{-h_{0}}}$ whose residues along $x_{h_{0}}, x_{-h_{0}}$ 
are $X^{\Delta_{2}}_{h_{0}}, X^{\Delta_{2}}_{-h_{0}}$ respectively by gluing 
$x_{h_{0}} \in P_{v_{h_{0}}}$ and $x_{-h_{0}} \in P_{v_{-h_{0}}}$. 
Therefore, by Theorem 3.1, 
$X^{\Delta_{1}}_{h}$ $(h \in \pm E \cup T)$ and $X^{\Delta_{2}}_{h'}$ $(h' \in \pm E' \cup T')$ 
satisfy the above rules. 
\ $\square$
\vspace{2ex}

\noindent
{\it 3.5. Unipotent periods of marked elliptic curves.} 
Let $\Delta = (V, E, T)$ be a stable graph of $(1, n)$-type, 
and assume that $\Delta$ is trivalent. 
Then by taking coordinates on $P_{v} = {\mathbb P}^{1}$ $(v \in V)$ 
such that these points corresponding to $h \in \pm E \cup T$ with $v_{h} = v$ 
belong to $\{ 0, 1, \infty \}$, 
$A_{\Delta} = {\mathbb Z}[[ y_{e} \, (e \in E)]]$. 
Furthermore, 
$\Delta$ can be obtained from the basic graph $\Delta_{0}$ 
by a combination of the alterations $\Delta_{1} \leftrightarrow \Delta_{2}$ in 2.2 
without shrinking the only one (oriented) loop $l$ in $\Delta_{0}$. 
Denote by $e_{l}$ the oriented edge in $\pm E$ obtained from $l$ under this operations, 
and put $E' = E - \{ |e_{l}| \}$. 
\vspace{2ex}

\noindent 
{\it Definition 3.3.} 
We define {\it unipotent periods} of ${\cal E}_{\Delta}$ as the monodromies of 
$\left( {\cal V}_{\Delta}, {\cal F}_{\Delta} \right)$ 
between ${\mathbb Z}$-rational tangential points 
on a degenerating family of $n$-marked elliptic curves over ${\mathbb C}$ 
obtained from ${\cal E}_{\Delta}^{\circ}$ as in 2.1 (P4). 
The unipotent periods are represented as noncommutative formal power series in 
$X^{\Delta}_{e}$ $(e \in E' \cup T)$. 
\vspace{2ex}

\noindent
{\bf Theorem 3.4.} 
\begin{it} 
Each coefficient of unipotent periods of ${\cal E}_{\Delta}$ 
as noncommutative formal power series in $X^{\Delta}_{e}$ $(e \in E' \cup T)$ 
is a formal power series in $y_{e}$  $(e \in E')$ whose coefficients are generated over ${\mathbb Q}$ 
by positive powers of $\pi \sqrt{-1}$, 
multiple zeta values and elliptic multiple zeta values. 
\end{it} 
\vspace{2ex}

\noindent
{\it Proof.} 
We will only prove the assertion in the case when $\Delta$ is obtained from a trivalent graph 
$\Delta_{1}$ with loop by a fusing move 
since one can prove the assertion in general cases by a similar method. 
Then $\Delta_{1} = (V_{1}, E_{1}, T_{1})$ has only one loop $l$ and nonloop edge $e_{1}$ 
with orientation such that $v_{\pm l} = v_{-e_{1}}$ and 
that $\Delta$ is obtained from $\Delta_{1}$ 
by the fusing move $e_{1} \rightarrow e'$ for certain edge $e'$ of $\Delta$. 

First, we consider unipotent periods of the family of $n$-marked complex curve of genus $1$  obtained from 
$$
\left( {\cal E}_{\Delta_{1}} \right)_{0} = 
{\cal E}_{\Delta_{1}} \otimes_{A_{\Delta_{1}}} A_{\Delta_{1}}/(y_{e_{1}}) 
$$ 
which is a union of the Tate curve ${\mathbb G}_{m} / \langle y_{l} \rangle$ 
and the universal rational curve $C_{0}$ with dual graph $\Delta_{1} - \{ l \}$. 
Then the unipotent periods for any rotation around 
${\mathbb Z}$-rational tangential points at $y_{e} = 0$ $(e \in E_{1})$ are 
integral powers of $\exp \left( \pi \sqrt{-1} X^{\Delta_{1}}_{e} \right)$. 
Furthermore, the cross ratios associated with the maximal degeneration of $C_{0}$ 
(cf. \cite[2.4]{BPP}) belong to ${\mathbb Z}[[ y_{e} \, (e \in E_{1} - \{ l, e_{1} \})]]$. 
Therefore, the assertion on $\left( {\cal E}_{\Delta_{1}} \right)_{0}$ follows from  \cite[Theorem 2.20]{BPP} and the definition of elliptic multiple zeta values 
reviewed in 3.3. 
  
Second, we consider unipotent periods of the family of $n$-marked elliptic curves 
over ${\mathbb C}$ obtained from ${\cal E}_{\Delta}$. 
This family is obtained by the fusing move $e_{1} \rightarrow e'$ 
corresponding to $0 \rightarrow 1 - s_{e'}$ in $(0, 1) \subset {\mathbb P}^{1}_{\mathbb R}$, 
where $s_{e'}$ can be regarded as a deformation parameter of ${\cal E}_{\Delta}$ 
associated with $e'$ by Theorem 2.2. 
Then this fusing move is represented as the monodromy along $(0, 1-s_{e'})$ 
of the connection 
$$
df - f \left( [T_{l}, A_{l}] \frac{dz}{z} + X_{e'} \frac{dz}{z-1} \right), 
$$
and hence the associated coefficients are seen in 3.2 to be formal power series in $s_{e'}$ 
whose coefficients are generated by multiple zeta vales over ${\mathbb Q}$. 
This fact together with the above assertion on $\left( {\cal E}_{\Delta_{1}} \right)_{0}$ 
imply the original assertion on ${\cal E}_{\Delta}$. 
\ $\square$

%\appendix
%\section{Some title}
%Please always give a title also for appendices.

%\acknowledgments

%This is the most common positions for acknowledgments. A macro is
%available to maintain the same layout and spelling of the heading.

%\paragraph{Note added.} This is also a good position for notes added
%after the paper has been written.


\renewcommand{\refname}{\normalsize{\bf References}}


\begin{thebibliography}{99}

\bibitem{BPP} 
Banks, P., Panzer, E., Pym, B.: 
Multiple zeta values in deformation quantization. 
to appear in Invent. Math. 
arXiv:1812.11649 

\bibitem{Br} 
Brown, F.: 
Multiple zeta values and periods of moduli spaces $\overline{\mathfrak M}_{0,n}$. 
Ann. Sci. \'{E}c. Norm. Sup\'{e}r. {\bf 42}, 371--489 (2009) 

%\bibitem{BoMW} 
%C. Bogner, S. M\"{u}ller-Stach and S. Weinzierl, 
%The unequal mass sunrise integral expressed through iterated integrals on 
%$\overline{\cal M}_{1,3}$, 
%Nuclear Phys. {\bf 954} (2020). 

\bibitem{CEE} 
Calaque, D., Enriquez, B., Etingof, P.: 
Universal KZB equations: the elliptic case. 
in Algebra, arithmetic, and geometry: in honor of Yu. I. Manin, vol. I, 
Progr. Math. vol. 269, Birkh\"{a}user, Boston, 2009, pp. 165--266 
arXiv: math/0702670 

\bibitem{D1} 
Deligne, P.: 
\'{E}quations diff\'{e}rentielles \`{a} points singuliers r\'{e}guliers. 
Lecture Notes in Math. vol. 163, Springer-Verlag, 1970 

\bibitem{D2} 
Deligne, P.: 
Le groupe fondamental de la droite projective moins trois points. 
in Galois groups over ${\mathbb Q}$, 
Publ. MSRI, vol. 16, Springer, 1989, pp. 79--298 

\bibitem{DM} 
Deligne, P., Mumford, D.: 
The irreducibility of the space of curves of given genus. 
Publ. Math. IHES {\bf 36}, 75--109 (1969) 

\bibitem{Dr} 
Drinfel'd, V. G.: 
On quasitriangular quasi-Hopf algebras and a group closely connected with 
${\rm Gal} \left( \overline{\mathbb Q}/{\mathbb Q} \right)$. 
Leningrad Math. J. {\bf 2}, 829--860 (1991) 

\bibitem{E} 
Enriquez, B.: 
Elliptic associators. 
Select. Math. {\bf 20}, 491--584 (2014) 
arXiv: 1003.1012

\bibitem{GD} 
Grothendieck, A., Dieudonn\'{e}, J.:  
El\'{e}ments de g\'{e}om\'{e}trie alg\'{e}brique: III, 
\'{E}tude cohomologique des faisceaux coh\'{e}rents, Premi\`{e}re partie. 
Publ. Math. IHES. {\bf 11} 5--167 (1961)

\bibitem{H} 
Hain, R.: 
Notes on the universal elliptic KZB equation. 
Pure and Applied Mathematics Quarterly, 
vol. 12, no. 2 (2016) International Press, 
arXiv: 1309.0580. 

\bibitem{I1} 
Ichikawa, T.: 
Generalized Tate curve and integral Teichm\"{u}ller modular forms, 
Amer. J. Math. {\bf 122} (2000) 1139--1174. 

%\bibitem{I2} 
%Ichikawa, T.: 
%Teichm\"{u}ller groupoids and Galois action. 
%J. reine angew. Math. {\bf 559}, 95--114 (2003) 

\bibitem{I3} 
Ichikawa, T.: 
The universal Mumford curves and its periods in arithmetic formal geometry. 
arXiv:2010.11517. 

\bibitem{IhN} 
Ihara, Y., Nakamura, H.: 
On deformation of maximally degenerate stable marked curves and Oda's problem. 
J. reine angew. Math. {\bf 487}, 125--151 (1997) 

\bibitem{K} 
Knudsen, F. F.: 
The projectivity of the moduli space of stable curves II, III. 
Math. Scand. {\bf 52}, 161--199, 200--212 (1983) 

\bibitem{KM} 
Knudsen, F. F., Mumford, D.: 
The projectivity of the moduli space of stable curves I. 
Math. Scand. {\bf 39}, 19--55 (1976) 

\bibitem{L} 
Luo, M.: 
The elliptic KZB connection and algebraic de Rham theory 
for unipotent fundamental groups of elliptic curves. 
Algebra Number Theory {\bf 13}, 2243--2275 (2019) 

%\bibitem{Mu} 
%Mumford, D.: 
%An analytic construction of degenerating curves over complete local rings. 
%Compos. Math. {\bf 24}, 129--174 (1972) 

%\bibitem{S} 
%Schottky, F.: 
%\"{U}ber eine specielle Function, welche bei einer bestimmten 
%inearen Transformation ihres Arguments unver\"{a}ndert bleibt. 
%J. reine angew. Math. {\bf 101}, 227--272 (1887) 

\bibitem{Si} 
Silverman, J. H.: 
Advanced topics in the arithmetic of elliptic curves. 
Graduate Texts in Math. vol. 151, Springer, 1994 

\bibitem{T} 
Tate, J.: 
A review of non-archimedean elliptic functions. 
in Elliptic Curves, Modular forms, \& Fermat's Last Theorem, 
International Press, Boston, 1995, pp. 162--184

\bibitem{W}
Weinzierl, S.: 
Iterated integrals related to Feynman integrals associated to elliptic curves, 
arXiv:2012.08429 

% Please avoid comments such as "For a review'', "For some examples",
% "and references therein" or move them in the text. In general,
% please leave only references in the bibliography and move all
% accessory text in footnotes.

% Also, please have only one work for each \bibitem.


\end{thebibliography}
\end{document}